\documentclass{amsart}
\usepackage{amssymb,amscd}

\newtheorem{lemma}{Lemma}
\newtheorem{prop}[lemma]{Proposition}
\newtheorem{cor}[lemma]{Corollary}
\newtheorem{thm}[lemma]{Theorem}

\newtheorem{thm?}[lemma]{Theorem?}

\DeclareFontEncoding{OT2}{}{} 
  \newcommand{\textcyr}[1]{%
    {\fontencoding{OT2}\fontfamily{wncyr}\fontseries{m}\fontshape{n}%
     \selectfont #1}}
\newcommand{\Sha}{{\mbox{\textcyr{Sh}}}}

\title{There are genus one curves of every index over every number field}
\author{Pete L. Clark}
\address{1126 Burnside Hall \\ Department of Mathematics and Statistics \\
McGill University \\ 805 Sherbrooke West \\ Montreal, QC, Canada H3A 2K6}
\email{clark@math.mcgill.ca}

\newcommand{\Q}{\ensuremath{\mathbb Q}}

\newcommand{\Z}{\ensuremath{\mathbb Z}}

\newcommand{\ra}{\ensuremath{\rightarrow}}

\newcommand{\Pic}{\operatorname{Pic}}

\newcommand{\Gm}{\mathbb G_m}

\newcommand{\Br}{\operatorname{Br}}

\newcommand{\Spec}{\operatorname{Spec}}

\newcommand{\KK}{\mathcal{K}}

\newcommand{\Jac}{\operatorname{Jac}}

\begin{document}
\maketitle

\begin{abstract}
We show that there exist genus one curves of every index over the rational numbers, answering affirmatively
a question of Lang and Tate.  The proof is ``elementary'' in the sense that it does not assume the finiteness of
any Shafarevich-Tate group.  On the other hand, using Kolyvagin's construction of a rational elliptic curve
whose Mordell-Weil and Shafarevich-Tate groups are both trivial, we show that there are infinitely many
curves of every index over every number field.
\end{abstract}
\section{Introduction}
\noindent
Let $C/K$ be a genus one curve over a field $K$.  There are two numerical invariants which quantify, in different ways,
the extent to which $C$ fails to have a $K$-rational point.  The \emph{index} of $C$ is the least degree of a field
extension $L/K$ such that $C$ has an $L$-rational point; equivalently, it is the least positive degree of a $K$-rational
divisor on $C$.  The \emph{period} of $C$ is the order of the cohomology class corresponding to $C$ in the Weil-Ch\^atelet
group $H^1(K,\Jac(C))$; equivalently, it is the least positive degree of a $K$-rational divisor class on $C$.  
It is well known (e.g. \cite{Lang-Tate}) that the period divides the index and that the two quantities have the same prime
divisors.
\\ \\
When $K$ is a number field, the index can strictly exceed the period \cite{Cassels}, \cite{Clark1}.  This is in a sense 
unfortunate, because while the index is of more geometric interest, it is the period which is directly addressed
by the machinery of Galois cohomology.  For example, it is an old theorem of Shafarevich-Cassels \cite[$\S$ 27]{Cassels2} that there are infinitely 
many classes of any given period in the Weil-Ch\^atelet group of any elliptic curve over a number field.  
\\ \indent With regard to the index, nearly half a century ago Lang and Tate 
\cite{Lang-Tate} asked the much more modest question of whether there are genus one curves of every index over $\Q$.  Lang and Tate were
able to show that if $E/K$ is an elliptic curve over a number field with a $K$-rational torsion point of order $n$, then $H^1(K,E)$ contains infinitely many
classes of index $n$.  In view of the uniform boundedness of torsion on elliptic curves, this does not get us very far.  Only very 
recently has 
substantial progress been made: in \cite{Stein}, W.A. Stein showed that for any number field $K$ there are infinitely many genus one curves over 
$K$ of index equal
to any number \emph{not divisible by $8$}.  
\\ \\
The following theorem and its corollary give a complete answer to the question of Lang and Tate.
\begin{thm}
Let $E/K$ be an elliptic curve over a number field with $E(K) = 0$.  For every positive integer
$n$, there exists an element $\eta \in H^1(K,E)$ of index $n$.
\end{thm}
\noindent
It is known (e.g. \cite{Nakagawa-Horie}) that there are infinitely many rational elliptic curves
with trivial Mordell-Weil group.  Thus:
\begin{cor}
For every positive integer $n$ there are infinitely many genus one curves $C/\Q$ of index $n$.
\end{cor}
\noindent
We want to emphasize that Theorem 1 does not require the finiteness of any Shafarevich-Tate group.  This is to 
be contrasted with the following results:
\begin{thm}
Let $E/K$ be an elliptic curve with $E(K) = 0$ and $\Sha(K,E) = 0$.  Then for every number field
$L/K$ and every positive integer $n$, there are infinitely many elements of $H^1(L,E)$ of index $n$.
\end{thm}
\noindent
Luckily for us, some examples of elliptic curves $E/\Q$ satisfying the hypotheses of Theorem 3 can be found
in a paper of Kolyvagin \cite[Theorem H]{Kolyvagin}: taking in his notation $D = -7$, we
get an elliptic curve $E/\Q$ with minimal Weierstrass equation $y^2 + y = x^3 - 49x - 86$ ($1813B1$ in
Cremona's tables) with $E(\Q) = \Sha(\Q,E) = 0$.  Thus we get the following, our main result:
\begin{cor}
There are infinitely many genus one curves of every index over every number field.
\end{cor}
\noindent
In Section 2 we set the stage with some preliminary results on a subset of $H^1(K,E)$ on which
the equality of period and index is guaranteed.  Readers familiar with the Heegner point
Euler system will recognize these classes as the (vastly simpler) analogue of classes constructed
by Kolyvagin.  In Section 3 we give the proofs of Theorems 1 and 3, and in Section 4 we make some
brief final remarks on generalizations and comparisons with Stein's work.
\section{The Kolyvagin set}
\noindent
For a number field $K$, we denote by $\Sigma_K$ the set of all places of $K$.
\\ \\
Let $E/K$ be an elliptic curve over a number field and $v \in \Sigma_K$.  We define 
$\mathcal{K}_v(K,E) \subset H^1(K,E)$ to be the subset of classes $\eta$ whose
local restriction to each $v' \neq v$ is zero.  We define the \emph{Kolyvagin set}
\[\mathcal{K}(K,E) = \bigcup_{v \in \Sigma_K} \mathcal{K}_v(K,E) \subset H^1(K,E). \]
The following proposition merely records for future reference some elementary properties of the Kolyvagin set.  The reader will have no
difficulty supplying the proof.
\begin{prop}
Let $E/K$ be an elliptic curve over a number field. \\
a) $\Sha(K,E) \subset \mathcal{K}(K,E)$. \\
b) For any $v$ and any positive integer $n$, $\KK_v(K,E)[n]$ is a finite group. \\
c) If $\eta \in \mathcal{K}(K,E)$ and $c \in \Z$, then $c \eta \in \mathcal{K}(K,E)$.  \\
\end{prop}
\noindent
The next result is the key observation about Kolyvagin classes that we will use to prove the main results of the paper.
\begin{prop}
Let $C/K$ be a genus one curve over a number field whose corresponding class
$\eta$ lies in $\KK(K,\Jac(C))$.  Then every rational divisor class on $C$ admits a rational divisor.  
In particular,
the period and index of $\eta$ are equal.
\end{prop}
\noindent
First proof: This follows rather formally from the existence of O'Neil's period-index obstruction map \cite{Cathy}.
Namely, for any elliptic curve $E/K$ defined over a field of characteristic zero and positive integer $n$, 
there exists a map \[\Delta: H^1(K,E[n]) \ra \Br(K)[n],\] functorial in $K$, and satisfying the following 
properties: a) a class $\eta \in H^1(K,E)[n]$ has index dividing $n$ if and only if there exists some
Kummer lift of $\eta$ to $\xi \in H^1(K,E[n])$ such that $\Delta(\xi) = 0$; b) if $\eta = 0$ then
every Kummer lift $\xi$ of $\eta$ has $\Delta(\xi) = 0$.  So let $\eta \in \KK(K,E)$ have period $n$, let 
$\xi$ be any lift of $\eta$ to $H^1(K,E[n])$ and consider $\Delta(\xi) \in \Br(K)[n]$.
Since for all $v' \neq v$, $\eta|_{v'} = 0$, by b) above we have that $\Delta(\xi)|_{v'} = 0$
in $\Br(K_{v'})$.  By virtue of the reciprocity law in the Brauer group of a number field,
we have that $\Delta(\xi) = 0$, so $\eta$ has index dividing $n$, hence index $n$.  \\ \\
Second proof: Let $V/K$ be any (smooth, projective geometrically irreducible) variety over any field 
$K$.  Taking low-degree terms in the Leray spectral sequence associated to the \'etale sheaf $\Gm$ on $\Spec K$ we
get \cite[Ch. IX]{BLR} an exact sequence
\[0 \ra \Pic(V/K) \ra \Pic(V/\overline{K})^{\mathfrak{g}_K} \stackrel{\delta}{\ra} \Br(K) \stackrel{\gamma}{\ra} \Br(V)  \]
which is functorial in $K$.  Thus the obstruction $\delta$ to a rational divisor class being represented
by a rational divisor is an element of the Brauer group of the base field $K$.  Moreover, a $K$-rational
point $P: \Spec K \ra V$ would induce a map $\Br(P) : \Br(V) \ra \Br(K)$ such that $\gamma \circ \Br(P) = 1_{\Br(V)}$, and it follows
that $V(K) \neq \emptyset$ implies $\delta \equiv 0$.  If now $K$ is a number field and $V$ is any variety which has 
rational points at every completion except possibly one, then the above reciprocity law argument gives us that $\delta \equiv 0$ on $V$.
\section{The proofs of Theorem 1 and Theorem 3}
\begin{lemma}
Let $E/K$ be an elliptic curve over a number field and $n$ be any positive integer.
If $v \in \Sigma_K$ is any finite place splitting completely in $K(E[n])$, then $H^1(K_v,E)$ has
an element of exact order $n$.
\end{lemma}
\noindent
Proof: By a seminal theorem of Tate \cite[Cor. I.3.4]{Milne}, the finite abelian groups $H^1(K_v,E)[n]$ and $E(K_v)/nE(K_v)$ are in duality, 
so it suffices to see that the latter group contains an element of exact order $n$ when $E$ has full $n$-torsion over $K_v$.
By the structure theory for compact $v$-adic Lie groups, $E(K_v) \cong \Z_{\ell}^N \oplus \Z/d_1\Z \oplus \Z/d_1d_2 \Z$
for some positive integers $N, \ d_1, \ d_2$; here $\ell$ is the residue characteristic of $K_v$.  If $E$ has
full $n$-torsion over $K_v$, then $n \ | \ d_1d_2$, so that any generator of $\Z/d_1d_2\Z$ has exact order
$n$ in $E(K)/nE(K)$.  
\\ \\
We now give the proof of Theorem 1, so let $E/K$ be an elliptic curve with $E(K) = 0$. \\
By primary decomposition for period and index of a cohomology class (e.g. \cite[Prop. 2.5]{Stein}),
it suffices to find classes of period and index equal to any prime power, say $n = p^a$.  There are
two cases to consider. \\ \\
Case 1: $\Sha(K,E)$ contains an element $\eta$ of exact order $p^a$.  Then by Propositions 5 and 6,
$\eta$ has index $p^a$.  \\ \\
Case 2: $\Sha(K,E)[p^{\infty}] = \Sha(K,E)[p^{a-1}]$ is a finite group.  By \cite[I.6.26(b)]{Milne},
whenever the $p$-primary torsion of the Shafarevich-Tate group of an abelian variety $A$ defined over a 
number field $K$ is a finite group (i.e., conjecturally always!) there is an exact sequence
\[0 \ra \Sha(K,A)[p^{\infty}] \ra H^1(K,A)[p^{\infty}] \ra \bigoplus_{v \in \Sigma_K} H^1(K_v,A)[p^{\infty}]
\ra (A^{\vee}(K)^{\wedge})^* \ra 0 , \]
where the three operations on the last term are, respectively, abelian variety dual, profinite completion,
and Pontrjagin dual.  But since we've assumed $E(K) = 0$, this gives a surjection 
\begin{equation}
H^1(K,E)[p^{\infty}] \ra \bigoplus_v H^1(K_v,E)[p^{\infty}] \ra 0.
\end{equation}
Invoking Lemma 7, let $v \in \Sigma_K$ be such that $H^1(K_v,E)$ contains an element
$\eta_{v}$ of exact order $p^{a}$.  
By (1), there
exists a global class $\eta$ which is locally trivial at every $v' \neq v$ and is locally equal to 
$\eta_{\ell}$, so that $\eta \in \KK(K,E)$.  By Proposition 4, any such $\eta$ has index equal to its period.
\\ \indent  The only
remaining question is what the period of $\eta$ is.  However, 
certainly the period of $\eta$ is of the form $c \cdot p^{a}$ for some positive integer $c$; then the class
$c \eta$ has exact period $p^{a}$. By Proposition 4c), $c \eta$ is still a Kolyvagin class, so also has
index $p^{a}$, as desired.  This completes the proof of Theorem 1.
\\ \\
We turn now the setting of Theorem 3, which is a sort of degenerate case of Theorem 1.  Indeed,
under the hypothesis that $\Sha(K,E) = 0$, the global duality relation becomes:
\[H^1(K,E) \stackrel{\sim}{\ra} \bigoplus_{v \in \Sigma_K} H^1(K_v,E). \]
It follows that the order of a Kolyvagin class is
always equal to the order of its nontrivial local restriction.  We finish, appropriately enough, by recalling
a result of Lang and Tate \cite[Cor. 1, p. 676]{Lang-Tate}: if $E/K_v$ is an elliptic curve over a $v$-adic field with good reduction and
$p^{a}$ is prime to the order of the residue field, then an element of $H^1(K_v,E)$ of period $p^{a}$ also has index 
$p^{a}$ and moreover is split by a local extension field $L_w/K_v$ if and only if $p^a$ divides the ramification index $e(L_w/K_v)$.  Thus if in Lemma 7 
we choose our place $v$ to be of good reduction for $E$, unramified in $L$ and prime to $p$ (which excludes only finitely many
places), we find that the class $\eta|_{L}$ has the property that $\eta_{L_w}$ has index at least $p^{a}$.
On the other hand since $\eta$ has index $p^{a}$, certainly $\eta|_{L}$ has index at most $p^{a}$, so
$\eta|_{L}$ has period and index both equal to $p^{a}$.  Since we can perform this construction for each $v$ in a set
of positive density, we get infinitely many such classes, completing the proof of Theorem 3.
\section{Some final thoughts}
\noindent
We should mention important work in progress of Stein and his students \cite{Stein2}, whose goal is to verify the full conjecture of 
Birch and Swinnerton-Dyer for all rational elliptic curves of conductor at most $1000$ and analytic rank at most one.  This work in particular
gives many other examples of rational elliptic curves satisfying the hypotheses of Theorem 3.  In fact there are close relations between
\cite{Stein}, \cite{Stein2} and the work of this paper, as we now explain.
\\ \\
The argument used to establish Theorem 3 readily gives the following variant:
\begin{thm}
Let $E/K$ be an elliptic curve over a number field whose Mordell-Weil group has rank zero.
Suppose moreover that $\Sha(K,E)$ is finite.  Then, for any $n$ prime to
$N = \#E(K) \cdot \# \Sha(K,E)$ and any number field $L/K$, there exist infinitely many
genus one curves $C/L$ of index $n$.
\end{thm}
\noindent
Theorem 8 should be compared with Stein's work \cite{Stein}, especially Theorem 3.1 and $\S 4.1$.  Stein works with the
rational elliptic curve $E = X_0(17)$.  The
calculations of \cite[$\S 5.1$]{Stein} showing that (in his notation) $B_K = 2$ give a good example of the techniques systemically employed
in \cite{Stein2}: namely, they show that $\Sha(\Q,X_0(17)[p] = 0$ for all odd primes $p$.  By a routine computation
with the $2$-Selmer group, one finds that $\Sha(\Q,X_0(17))[2] = 0$.  Thus $\Sha(\Q,X_0(17)) = 0$, as predicted by BSD.  However $X_0(17)(\Q) = 
\Z/4\Z$, so by Theorem 8, for all number fields $K$ there are elements of $H^1(K,X_0(17))$ of every \emph{odd} index.
\\ \\
Much more remains to be done on the index problem for genus one curves over number fields.  The obvious analogue for the index of 
the Shafarevich-Cassels theorem is the assertion that for any positive integer $n$ and any elliptic curve $E$ over a number field $K$, there exist 
infinitely many classes in $H^1(K,E)$ with index $n$.  The case of $n=2$ (``biconic curves'') and arbitrary $E$ and $K$ can be handled using the methods of 
\cite{Clark1} together with the theory of explicit $2$-descent; this is the subject of work in progress of the present author.  The general case seems to be 
quite difficult and will probably require new ideas.  

\end{document}